\documentclass[a4paper,12pt,twoside]{article}
\usepackage{times}
\usepackage[english]{babel}
\usepackage{amssymb,amsmath}
\numberwithin{equation}{section}
\usepackage{amsthm}
\usepackage{array}
\usepackage[noadjust]{cite}
\usepackage{hyperref}
\usepackage[height=22.7cm , width = 16cm , top = 4cm , left = 3cm, a4paper]{geometry}
\usepackage{amssymb}
\usepackage{amsmath}
\usepackage{cite}
 \usepackage{pstricks}
\usepackage{epsfig}
\usepackage{verbatim}
\usepackage{multicol}
\usepackage{graphicx, color, psfrag}
\usepackage{graphics, psfrag}
\usepackage{color,soul}
\newtheorem{theorem}{Theorem}[section]
\newtheorem{remark}[theorem]{Remark}
\newcommand{\e}{\end{document}}

\theoremstyle{definition}
\newtheorem{definition}{Definition}[section]
\begin{document}

\thispagestyle{empty}

\author{F. M. Abdel Moneim$^{(1)}$\thanks{ Corresponding author, email: \href{mailto:fatmakasem1982@gmail.com}{fatmakasem1982@gmail.com}},~ Abdelfattah Mustafa $^{(1,2)}$\thanks{Email: \href{mailto: amelsayed@mans.edu.eg}{amelsayed@mans.edu.eg}},~ and B. S. El-Desouky$^{(1)}$ \\
$^{(1)}$ {\small Department of Mathematics, Faculty of Science,
Mansoura University, Mansoura 35516, Egypt.} \\
$^{(2)}$ {\small Department of Mathematics, Faculty of Science, Islamic University of Madinah, KSA.}}

\title{New Generalization Families of Higher Order Daehee  Numbers and Polynomials}

\date{}

\maketitle
\small \pagestyle{myheadings}
        \markboth{{\scriptsize New Generalization Families of Higher Order Daehee Numbers and Polynomials}}
        {{\scriptsize {F. M. Abdel Moneim, A. Mustafa and B. S. El-Desouky}}}

\hrule \vskip 8pt

\begin{abstract}
In this paper, we present a new definition  and generalization of higher order Daehee of the first and second kind. Some new results for these polynomials and numbers  are derived. Furthermore, some interesting special cases of the  new generalized higher order Daehee numbers and polynomials are deduced.

\end{abstract}

\noindent
{\bf Keywords:}
{\it } Daehee  numbers; Daehee  polynomials; higher-order Daehee  numbers; higher-order Daehee polynomials.

\noindent
{\it {\bf AMS Subject Classification:}} {\em   05A19; 11B73; 11T06.}

\section{Introduction}
The $n$-th Daehee polynomials are defined by the generating function, \cite{10,11,12,13,14,15,16,17,18},
\begin{equation} \label{1}
\left( \frac{\log(1+t)}{t} \right) (1+t)^x = \sum_{n=0}^\infty D_n (x)  \frac{t^n}{n!} .
\end{equation}

\noindent
In the special case, $x=0,D_n=D_n(0)$ are called Daehee numbers, and for $n\geq 0$,
\begin{equation} \label{2}
\int_{\mathbb{Z}_p} (x )_n  d\mu_0 (x)=D_n.
\end{equation}

\noindent
For $k\in \mathbb{N}$,  Kim \cite{10} introduced Daehee numbers of the first kind of order $k$ by
\begin{equation} \label{3}
D_n^{(k) }=\int_{\mathbb{Z}_p} \int_{\mathbb{Z}_p} \cdots \int_{\mathbb{Z}_p} (x_1+x_2+ \cdots x_k )_n  d\mu_0 (x_1 ) d\mu_0 (x_2 )\cdots d\mu_0 (x_k ),
\end{equation}

\noindent
where $n$ is nonnegative integer.\\

\noindent
The generating function of these numbers are given as
\begin{equation}
\sum_{n=0}^\infty D_n^{(k)}  \frac{t^n}{n!}  =\left(\frac{\log(1+t)}{t}\right)^k,
\end{equation}

\noindent
where $n\in Z \geq 0, k\in N$.\\

\noindent
The higher-order Daehee polynomials are defined as
\begin{equation} \label{5}
D_n^{(k)} (x)=\int_{\mathbb{Z}_p}  \int_{\mathbb{Z}_p} \cdots  \int_{\mathbb{Z}_p} (x_1+x_2+ \cdots x_k+x)_n  d\mu_0 (x_1 )d\mu_0 (x_2 )\cdots d\mu_0 (x_k ).
\end{equation}

\noindent
For $k\in \mathbb{Z}$, the Bernoulli polynomials of order $k$ are defined by the generating function to be, \cite{1,2,3,10},
\begin{equation} \label{6}
\left(\frac{t}{e^t-1} \right)^k e^{xt}= \sum_{n=0}^\infty B_n^{(k)}(x) \frac{t^n}{n!},
\end{equation}

\noindent
where $x=0,B_n^{(k) }=B_n^{(k)} (0)$ are called the Bernoulli numbers of order $k$. Also, Kim proved that

\begin{equation} \label{7}
D_n^{(k)} (x)= \sum_{l=0}^n s(n,\ell)  B_\ell^{(k)} (x),
\end{equation}
and
\begin{equation} \label{8}
B_n^{(k)} (x)= \sum_{\ell=0}^n S(n,\ell) D_\ell^{(k) } (x).
\end{equation}

\noindent
An explicit formula for higher-order Daehee numbers are given by
\begin{equation} \label{9}
D_n^{(k)}=\frac{s(n+k,k)}{\binom{n+k}{k}}, \quad n\geq  0, k\geq 1,
\end{equation}

\noindent
where $s(n+k,k)$ are Stirling numbers of the first kind see \cite{8,10}.

\noindent
Comtet \cite{2} introduced two kinds of Cauchy numbers, the first kind of Cauchy number is given by
\begin{equation} \label{10}
C_n=\int_0^1(x )_n  dx,
\end{equation}

\noindent
and the second kind is given by
\begin{equation} \label{11}
\hat{C}_n=\int_0^1(-x )_n  dx.
\end{equation}

\section{Generalized Higher Order Daehee Numbers and Polynomials of the First Kind}
In this section, we introduce new definitions of the generalization families for Daehee numbers and polynomials. Also, some new results for these numbers and polynomials are derived. Some special cases are discussed.

\begin{definition} \label{def:2.1}
The generalized higher order Daehee numbers  of the first kind, $\check{D}_{n;\bar{\pmb \alpha}, \bar{\pmb r}}^{(k)}$,  are defined by
\begin{equation} \label{12}
\check{D}_{n;\bar{\pmb \alpha}, \bar{\pmb r}}^{(k)}=\int_{\mathbb{Z}_p} \int_{\mathbb{Z}_p} \cdots  \int_{\mathbb{Z}_p} \prod_{i=0}^{n-1} (x_1 x_2\cdots x_k-\alpha_i)^{r_i }   d\mu_0 (x_1 )d\mu_0 (x_2 )\cdots d\mu_0 (x_k ).
\end{equation}
\end{definition}

\noindent
Some new results for  $\check{D}_{n;\bar{\pmb \alpha}, \bar{\pmb r}}^{(k)}$ can be derived as follows.

\begin{theorem} \label{Th:2.1}
The numbers $\check{D}_{n;\bar{\pmb \alpha}, \bar{\pmb r}}^{(k)}$ satisfy the relation
\begin{equation} \label{13}
\check{D}_{n;\bar{\pmb \alpha}, \bar{\pmb r}}^{(k)}= \sum_{m=0}^{|r|}   s_{\bar{\pmb \alpha}}(n,m; \bar{\pmb r})  \sum_{\ell_1=0}^m \cdots \sum_{\ell_k=0}^m \prod_{i=0}^k S(m,\ell_i )  D_{\ell_i}.
\end{equation}

\end{theorem}

\begin{proof}
The generalized Comtet numbers of the first and second kind, where $s_{\bar{\pmb \alpha}}  (n,i;\bar{\pmb r} )$,  $S_{\bar{\pmb \alpha}}  (n,i;\bar{\pmb r})$, see \cite{2}, are defined, respectively, by

\begin{equation} \label{14}
(x;\bar{\pmb \alpha},\bar{\pmb r})_n=\sum_{i=0}^n  s_{\bar{\pmb \alpha}}  (n,i;\bar{\pmb r})   x^i,
\end{equation}
and
\begin{equation} \label{15}
x^n= \sum_{i=0}^n  S_{\bar{\pmb \alpha}} (n,i;\bar{\pmb r}) (x;\bar{\pmb \alpha},\bar{\pmb r})_i,
\end{equation}

\noindent
where $(x;\bar{\pmb \alpha},\bar{\pmb r})_n=\prod_{i=0}^{n-1}(x-\alpha_i  )^{r_i}$, $\bar{\pmb \alpha}=(\alpha_0,\alpha_1,\cdots,\alpha_{n-1} )$ and $\bar{\pmb r}=(r_0,r_1,\cdots,r_{n-1} )$.\\

\noindent
From Eq. \eqref{12}, we have
\begin{eqnarray} \label{16}
\check{D}_{n;\bar{\pmb \alpha},\bar{\pmb r}}^{(k)} &= &
\int_{\mathbb{Z}_p}  \int_{\mathbb{Z}_p} \cdots \int_{\mathbb{Z}_p} \sum_{m=0}^{|r|} s_{\bar{\pmb \alpha}}  (n,m;\bar{\pmb r} )  (x_1 x_2\cdots x_k  )^m  d\mu_0 (x_1 )\cdots d\mu_0 (x_k )
\\
&= & \sum_{m=0}^{|r|} s_{\bar{\pmb \alpha}}  (n,m;\bar{\pmb r})  \int_{\mathbb{Z}_p}  \int_{\mathbb{Z}_p}\cdots \int_{\mathbb{Z}_p} (x_1 x_2\cdots x_k  )^m  d\mu_0 (x_1 )\cdots d\mu_0 (x_k )
\nonumber\\
&= & \sum_{m=0}^{|r|} s_{\bar{\pmb \alpha}}  (n,m;\bar{\pmb r}) \int_{\mathbb{Z}_p} (x_1  )^m  d\mu_0 (x_1 ) \cdots \int_{\mathbb{Z}_p} (x_k  )^m  d\mu_0 (x_k )
\nonumber\\
& =& \sum_{m=0}^{|r|} s_{\bar{\pmb \alpha}}  (n,m;\bar{\pmb r})  \left[ \int_{\mathbb{Z}_p} \sum_{\ell_1=0}^m S(m,\ell_1 )  (x_1  )_{\ell_1}  d\mu_0 (x_1 )\cdots \int_{\mathbb{Z}_p} \sum_{\ell_k=0}^m S(m,\ell_k )  (x_k )_{\ell_k}  d\mu_0 (x_k )\right]
\nonumber\\
&=& \sum_{m=0}^{|r|} s_{\bar{\pmb \alpha}}  (n,m;\bar{\pmb r})  \left[ \sum_{\ell_1=0}^ S(m,\ell_1 )  \int_{\mathbb{Z}_p} (x_1 )_{\ell_1}  d\mu_0 (x_1 )\cdots \sum_{\ell_k=0}^m S(m,\ell_k )  \int_{\mathbb{Z}_p} (x_k )_{\ell_k}  d\mu_0 (x_k )\right]
\nonumber\\
&=& \sum_{m=0}^{|r|} s_{\bar{\pmb \alpha}}  (n,m;\bar{\pmb r})  \left[ \sum_{\ell_1=0}^m S(m,\ell_1 )  D_{\ell_1} \cdots \sum_{\ell_k=0}^m S(m,\ell_k )  D_{\ell_k} \right]
\nonumber\\
&=& \sum_{m=0}^{|r|}  s_{\bar{\pmb \alpha}}  (n,m;\bar{\pmb r})  \left[\sum_{\ell_1=0}^m \sum_{\ell_2=0}^m \cdots \sum_{\ell_k=0}^m S(m,\ell_1 ) S(m,\ell_2 )\cdots S(m,\ell_k )  D_{\ell_1 } \cdots D_{\ell_k } \right].
\nonumber
\end{eqnarray}

\noindent
Then, we obtain \eqref{13}.
\end{proof}

\begin{theorem} \label{Th:2.2}
The numbers $\check{D}_{n;\bar{\pmb \alpha}, \bar{\pmb r}}^{(k)}$  satisfy the relation
\begin{equation} \label{17}
\check{D}_{n;\bar{\pmb \alpha},\bar{\pmb r}}^{(k) }= \sum_{m=0}^{|r|}   s_{\bar{\pmb \alpha}} (n,m;\bar{\pmb r})  \sum_{\ell_1=0}^m \cdots \sum_{\ell_k=0}^m \prod_{i=0}^k \frac{(-1)^{\ell_i } \ell_i! S(m,\ell_i )}{\ell_i+1},
\end{equation}

\noindent
which gives a relationship of the generalized higher order Daehee numbers of the first kind in terms of the multiparameter non-central Stirling numbers of the second kind and Stirling number of the first kind, see \cite{4,8,9,10}.
\end{theorem}

\begin{proof}
Substituting Eq. \eqref{9} into \eqref{13}, we obtain \eqref{17}.
\end{proof}

\begin{remark} \label{Rem:1}
\begin{eqnarray} \label{18}
\int_{\mathbb{Z}_p} \int_{\mathbb{Z}_p} \cdots  \int_{\mathbb{Z}_p} \prod_{i=0}^{n-1} (x_1 x_2\cdots x_k-\alpha_i  )^{r_i }  d\mu_0 (x_1 )d\mu_0 (x_2 )\cdots d\mu_0 (x_k )
\nonumber\\
=\sum_{m=0}^{|r|}   s_{\bar{\pmb \alpha}}  (n,m;\bar{\pmb r}) \sum_{\ell_1=0}^m \cdots \sum_{\ell_k=0}^m \prod_{i=0}^k \frac{(-1)^{\ell_i} \ell_i ! S(m,\ell_i )}{\ell_i+1}.
\end{eqnarray}
\end{remark}

\begin{definition} \label{def:2.2}
The multi-parameter poly-Cauchy numbers of the first kind  $C_{n;\bar{\pmb \alpha},\bar{\pmb r}}^{(k)}$ are defined by, \cite{6}
\begin{equation} \label{19}
C_{n;\bar{\pmb \alpha},\bar{\pmb r}}^{(k) }=\int_0^{\ell_1} \int_0^{\ell_2} \cdots \int_0^{\ell_k} \prod_{i=0}^{n-1} (x_1 x_2 \cdots x_k-\alpha_i  )^{r_i}  dx_1 dx_2\cdots dx_k.
\end{equation}
\end{definition}

\begin{theorem} \label{Th:2.3}
The numbers $C_{n;\bar{\pmb \alpha},\bar{\pmb r}}^{(k)}$   satisfy the  relation
\begin{equation} \label{20}
C_{n;\bar{\pmb \alpha},\bar{\pmb r}}^{(k)}=\sum_{i=0}^{|r|}  s_{\bar{\pmb \alpha}}  (n,i;\bar{\pmb r} )   \frac{ (\ell_1 \ell_2 \cdots \ell_k )^{j+1}}{(j+1)^k}.
\end{equation}
\end{theorem}

\begin{proof}
From Eq. \eqref{19}, we have
\begin{eqnarray} \label{21}
C_{n;\bar{\pmb \alpha},\bar{\pmb r}}^{(k)}& = & \int_0^{\ell_1}  \int_0^{\ell_2} \cdots \int_0^{\ell_k} \sum_{m=0}^{|r|}  s_{\bar{\pmb \alpha}}  (n,m;\bar{\pmb r})  (x_1 x_2 \cdots x_k  )^m  dx_1 dx_2 \cdots dx_k
\nonumber\\
&=& \sum_{m=0}^{|r|}  s_{\bar{\pmb \alpha}} (n,m;\bar{\pmb r}) \int_0^{\ell_1} \int_0^{\ell_2} \cdots \int_0^{\ell_k} (x_1 x_2 \cdots x_k  )^m  dx_1 dx_2 \cdots dx_k
\nonumber\\
&=& \sum_{m=0}^{|r|}s_{\bar{\pmb \alpha}}  (n,m;\bar{\pmb r} ) \int_0^{\ell_1} (x_1  )^m dx_1  \int_0^{\ell_1} (x_1  )^m dx_1  \cdots \int_0^{\ell_k}(x_k  )^m dx_k,
\end{eqnarray}

\noindent
then, we obtain \eqref{20}.
\end{proof}

\begin{definition} \label{Def:2.3}
The generalized Daehee numbers $\check{D}_n^{(k)}$  are defined by
\begin{equation} \label{22}
\check{D}_n^{(k)} =\int_{\mathbb{Z}_p}  \int_{\mathbb{Z}_p} \cdots \int_{\mathbb{Z}_p} (x_1 x_2 \cdots x_k  )_n  d\mu_0 (x_1 )d\mu_0 (x_2 ) \cdots d\mu_0 (x_k ).
\end{equation}
\end{definition}

\begin{theorem} \label{Th:2.4}
The polynomials $\check{D}_{n;\bar{\pmb \alpha},\bar{\pmb r}}^{(k)}$ satisfy the relation
\begin{equation} \label{23}
\check{D}_{n;\bar{\pmb \alpha},\bar{\pmb r}}^{(k)} (x)= \sum_{\ell=0}^{|r|} S(n,\ell;\bar{\pmb \alpha},\bar{\pmb r})   \check{D}_\ell^{(k)}.
\end{equation}
\end{theorem}

\begin{proof}
From Eq. \eqref{21}, we have

\begin{eqnarray}\label{24}
\check{D}_{n;\bar{\pmb \alpha},\bar{\pmb r}}^{(k)} (x) & = &
\int_{\mathbb{Z}_p} \cdots \int_{\mathbb{Z}_p} \sum_{m=0}^{|r|} s_{\bar{\pmb \alpha}} (n,m;\bar{\pmb r}) (x_1 x_2\cdots x_k x)^m  d\mu_0 (x_1) \cdots d\mu_0 (x_k )
\nonumber\\
&=& \sum_{m=0}^{|r|} s_{\bar{\pmb \alpha}}  (n,m;\bar{\pmb r}) \int_{\mathbb{Z}_p} \int_{\mathbb{Z}_p} \cdots \int_{\mathbb{Z}_p} (x_1 x_2 \cdots x_k  )^m  d\mu_0 (x_1) \cdots d\mu_0 (x_k )
\nonumber\\
&= & \sum_{m=0}^{|r|} s_{\bar{\pmb \alpha}}  (n,m;\bar{\pmb r} )  \int_{\mathbb{Z}_p} \int_{\mathbb{Z}_p} \cdots \int_{\mathbb{Z}_p} \sum_{\ell=0}^m S(m,\ell)  (x_1 x_2 \cdots x_k )_\ell  d\mu_0 (x_1) \cdots d\mu_0 (x_k )
\nonumber\\
&=&\sum_{m=0}^{|r|} s_{\bar{\pmb \alpha}}  (n,m;\bar{\pmb r}) \sum_{\ell=0}^m S(m,\ell)  \int_{\mathbb{Z}_p}  \int_{\mathbb{Z}_p} \cdots \int_{\mathbb{Z}_p} (x_1 x_2 \cdots x_k)_\ell  d\mu_0 (x_1 )\cdots d\mu_0 (x_k ).
\end{eqnarray}

\noindent
Substituting from Eq. \eqref{22} into Eq. \eqref{24}, we have
\begin{equation*}
\check{D}(n;\bar{\pmb \alpha},\bar{\pmb r})^{(k)} (x)=\sum_{m=0}^{|r|} s_{\bar{\pmb \alpha}}  (n,m;\bar{\pmb r} ) \sum_{\ell=0}^m S(m,\ell)  \check{D}_\ell^{(k) }=\sum_{\ell=0}^{|r|} \sum_{m=\ell}^{|r|} s_{\bar{\pmb \alpha}}  (n,m;\bar{\pmb r} )S(m,\ell) \check{D}_\ell^{(k) },
\end{equation*}

\noindent
since,
\begin{equation*}
\sum_{m=\ell}^{|r|} s_{\bar{\pmb \alpha}}(n,m;\bar{\pmb r})S(m,\ell) =S(n,\ell;\bar{\pmb \alpha},\bar{\pmb r} ),
\end{equation*}

\noindent
where $S(n,\ell;\bar{\pmb \alpha},\bar{\pmb r})$ are the generalized multiparameter non central Stirling numbers of the second kind, see \cite[Eq. (4.4)]{2}, hence, we obtain Eq. \eqref{23}.
\end{proof}


\noindent
The  higher order Daehee polynomials of the first kind,  $\check{D}_{n;\bar{\pmb \alpha },\bar{\pmb r }}^{(k) }(x)$, can be defined as follows.

\begin{definition} \label{Def:3.1}
The generalized higher order Daehee polynomials of the first kind,  $\check{D}_{n;\bar{\pmb \alpha },\bar{\pmb r }}^{(k) }(x)$ are defined by
\begin{equation} \label{25}
\check{D}_{n;\bar{\pmb \alpha },\bar{\pmb r }}^{(k)} (x)=\int_{\mathbb{Z}_p} \int_{\mathbb{Z}_p} \cdots \int_{\mathbb{Z}_p} \prod_{i=0}^{n-1} (x_1 x_2 \cdots x_k x-\alpha_i  )^{r_i }  d\mu_0 (x_1 )d\mu_0 (x_2 )\cdots d\mu_0 (x_k ).
\end{equation}
\end{definition}

\begin{definition} \label{3.2}
The generalized Daehee polynomials $\check{D}_n^{(k)}(x)$ are defined by
\begin{equation} \label{26}
\check{D}_n^{(k)} (x)=\int_{\mathbb{Z}_p} \int_{\mathbb{Z}_p} \cdots \int_{\mathbb{Z}_p} (x_1 x_2 \cdots x_k x )_n  d\mu_0 (x_1 ) d\mu_0 (x_2 ) \cdots d\mu_0 (x_k ).
\end{equation}
\end{definition}

\begin{theorem} \label{Th:3.1}
The polynomials $\check{D}_{n; \bar{\pmb \alpha },\bar{\pmb r }}^{(k) } (x)$ satisfy the relation
\begin{equation} \label{27}
\check{D}_{n; \bar{\pmb \alpha }, \bar{\pmb r }}^{(k) } (x)= \sum_{\ell=0}^{|r|} S(n,\ell; \bar{\pmb \alpha }, \bar{\pmb r } )   \check{D}_\ell^{(k)}  (x).
\end{equation}
\end{theorem}

\begin{proof}
From Eq. \eqref{25}, we have
\begin{equation*}
\check{D}_{n; \bar{\pmb \alpha }, \bar{\pmb r }}^{(k)} (x) = \int_{\mathbb{Z}_p} \cdots  \int_{\mathbb{Z}_p}  \sum_{m=0}^{|r|} s_{\bar{\pmb \alpha }}(n,m; \bar{\pmb r } )  (x_1 x_2 \cdots x_k x )^m  d\mu_0 (x_1 ) \cdots d\mu_0 (x_k ).
\end{equation*}

\noindent
Using Theorem 2.4, hence, we obtain Eq. \eqref{27}.
\end{proof}

\subsection{Some special cases}
The Daehee numbers and polynomials of the first kind, \cite{9,10,13}, can be obtained from the new definition as a special cases.

\noindent
{\bf Case 1:} (i) Setting  $r_i=r,\alpha_i=i$  in Eq. \eqref{25}, we have
\begin{eqnarray}
\check{D}_{n;i,r}^{(k)} (x) &=&
 \int_{\mathbb{Z}_p} \int_{\mathbb{Z}_p} \cdots \int_{\mathbb{Z}_p} \prod_{i=0}^{n-1} (x_1 x_2 \cdots x_k x-i )^r  d\mu_0 (x_1 )d\mu_0 (x_2 )\cdots d\mu_0 (x_k )
\nonumber\\
&=& \int_{\mathbb{Z}_p}  \int_{\mathbb{Z}_p} \cdots \int_{\mathbb{Z}_p} (x_1 x_2 \cdots x_k x )_{nr} \; d\mu_0 (x_1 )d\mu_0 (x_2 )\cdots d\mu_0 (x_k).
\end{eqnarray}

\noindent
Replacing  $nr$ by $n$, we obtain the generalized Daehee polynomials.\\

\noindent
(ii) Setting  $r_i=r,\alpha_i=i$  in Eq. \eqref{12}, we obtain
\begin{equation} \label{28}
\check{D}_{n;i,r }^{(k) }=\int_{\mathbb{Z}_p} \int_{\mathbb{Z}_p} \cdots \int_{\mathbb{Z}_p} (x_1 x_2\cdots x_k  )_{nr}\;   d\mu_0 (x_1 ) d\mu_0 (x_2 )\cdots d\mu_0 (x_k ).
\end{equation}

\noindent
Replacing  $nr$ by $n$, we obtain the higher-order Daehee numbers.\\

\noindent
{\bf Case 2:} (i) Setting  $r_i=r,\alpha_i=\alpha$ in Eq. \eqref{25}, we obtain
\begin{eqnarray} \label{29}
\check{D}_{n;\alpha,r}^{(k)} (x)&=&
\int_{\mathbb{Z}_p}   \int_{\mathbb{Z}_p} \cdots \int_{\mathbb{Z}_p} (x_1 x_2 \cdots x_k x-\alpha )^{nr}  d\mu_0 (x_1 )d\mu_0 (x_2 )\cdots d\mu_0 (x_k )
\nonumber\\
&= & \int_{\mathbb{Z}_p} \int_{\mathbb{Z}_p} \cdots \int_{\mathbb{Z}_p} \sum_{\ell=0}^{nr} S(nr,\ell)  (x_1 x_2 \cdots x_k x-\alpha )_\ell  d\mu_0 (x_1 )d\mu_0 (x_2 )\cdots d\mu_0 (x_k )
\nonumber\\
&= & \sum_{\ell=0}^{nr} S(nr,\ell)  \int_{\mathbb{Z}_p}  \int_{\mathbb{Z}_p} \cdots \int_{\mathbb{Z}_p}(x_1 x_2 \cdots x_k x-\alpha )_\ell d\mu_0 (x_1 )d\mu_0 (x_2 )\cdots d\mu_0 (x_k )
\nonumber\\
&= & \sum_{\ell=0}^{nr} S(nr,\ell) \check{D}_{\ell,\alpha }^{(k) } (x).
\end{eqnarray}

\noindent
(ii) Setting  $r_i=r,\alpha_i=\alpha$ in Eq. \eqref{12}, we obtain
\begin{eqnarray} \label{30}
\check{D}_{n;\alpha,r}^{(k)} &=&
\int_{\mathbb{Z}_p} \int_{\mathbb{Z}_p} \cdots  \int_{\mathbb{Z}_p} (x_1 x_2 \cdots x_k-\alpha )^{nr}  d\mu_0 (x_1 ) d\mu_0 (x_2 ) \cdots d\mu_0 (x_k )
\nonumber\\
&= & \sum_{\ell=0}^{nr} S(nr,\ell) \check{D}_{\ell,\alpha}^{(k) }.
\end{eqnarray}

\noindent
At $\alpha_i=0$ in Eq. \eqref{29}, we obtain
\begin{equation} \label{31}
\check{D}_{n;0,r}^((k) ) (x)=\sum_{\ell=0}^{nr} S(nr,\ell) \check{D}_\ell^{(k) } (x).
\end{equation}

\noindent
At $\alpha_i=0$ in Eq. \eqref{31}, we obtain
\begin{equation} \label{32}
\check{D}_{n;0,r}^{(k) }=\sum_{\ell=0}^{nr} S(nr,\ell) \check{D}_\ell^{(k)}.
\end{equation}

\noindent
{\bf Case 3:} (i) Setting  $r_i=1,\alpha_i=\alpha$ in Eq. \eqref{25}, we obtain
\begin{equation} \label{33}
\check{D}_{n;\alpha,1}^{(k)} (x)=\sum_{\ell=0}^n S(n,\ell) \check{D}_{\ell,\alpha}^{(k)} (x).
\end{equation}

\noindent
(ii) Setting  $r_i=1,\alpha_i=\alpha$ in Eq. \eqref{12}, we obtain

\begin{equation} \label{34}
\check{D}_{n;\alpha,1}^{(k)}=\sum_{\ell=0}^n S(n,\ell) \check{D}_{\ell,\alpha}^{(k) }.
\end{equation}

\noindent
(iii) Setting  $r_i=1,  \alpha_i=1$ in Eq. \eqref{25}, we obtain
\begin{eqnarray} \label{35}
\check{D}_{n;1,1}^{(k)} (x)& = &
\int_{\mathbb{Z}_p} \cdots  \int_{\mathbb{Z}_p} (x_1 x_2 \cdots x_k x-1 )^n  d\mu_0 (x_1 )d\mu_0 (x_2 )\cdots d\mu_0 (x_k )
\nonumber\\
&=& \sum_{\ell=0}^n S(n,\ell) \check{D}_{\ell,1}^{(k)} (x).
\end{eqnarray}

\noindent
(iv) Setting  $r_i=1,   \alpha_i=1$ in Eq. \eqref{12}, we obtain
\begin{eqnarray} \label{36}
\check{D}_{n;1,1}^{(k) } &=&
\int_{\mathbb{Z}_p} \cdots \int_{\mathbb{Z}_p} (x_1 x_2 \cdots x_k-1 )^n  d\mu_0 (x_1 )d\mu_0 (x_2 )\cdots d\mu_0 (x_k )
\nonumber\\
&=&
\sum_{\ell=0}^n S(n,\ell) \check{D}_{\ell,1}^{(k) }.
\end{eqnarray}

\noindent
{\bf Case 4:} (i) Setting  $r_i=1,\alpha_i=0$ in Eq. \eqref{25}, we obtain
\begin{eqnarray} \label{37}
\check{D}_{n;0,1}^{(k)} (x)& =&
\int_{\mathbb{Z}_p} \int_{\mathbb{Z}_p} \cdots \int_{\mathbb{Z}_p} (x_1 x_2 \cdots x_k x)^n  d\mu_0 (x_1 )d\mu_0 (x_2 )\cdots d\mu_0 (x_k )
\nonumber\\
&=& \sum_{\ell=0}^n S(n,\ell) \check{D}_\ell^{(k)} (x).
\end{eqnarray}

\noindent
(ii) Setting  $r_i=1,\alpha_i=0$ in Eq. \eqref{12}, we obtain
\begin{equation} \label{38}
\check{D}_{n;0,1}^{(k)}=\sum_{\ell=0}^n S(n,\ell) \check{D}_\ell^{(k) }.
\end{equation}

\noindent
{\bf Case 5:} (i) Setting  $r_i=1, \alpha_i=i,\; i=0,1,\cdots, n-1$  in Eq. \eqref{25}, we have
\begin{eqnarray} \label{39}
\check{D}_{n;i,1}^{(k)} (x) & = &
\int_{\mathbb{Z}_p} \int_{\mathbb{Z}_p} \cdots \int_{\mathbb{Z}_p} \prod_{i=0}^{n-1} (x_1 x_2 \cdots x_k x-i )  d\mu_0 (x_1 )d\mu_0 (x_2 )\cdots d\mu_0 (x_k )
\nonumber\\
&= & \int_{\mathbb{Z}_p}  \int_{\mathbb{Z}_p} \cdots \int_{\mathbb{Z}_p} (x_1 x_2 \cdots x_k x)_n  d\mu_0 (x_1 )d\mu_0 (x_2 ) \cdots d\mu_0 (x_k )=\check{D}_n^{(k) } (x),
\end{eqnarray}

\noindent
we obtain the higher-order Daehee polynomials which defined by Kim, see \cite{10}.\\

\noindent
(ii) Setting  $r_i=1, \alpha_i=i,\; i=0,1,\cdots, n-1$  in Eq. \eqref{12}, we have
\begin{eqnarray} \label{40}
\check{D}_{n;i,1}^{(k)} &= &
\int_{\mathbb{Z}_p} \int_{\mathbb{Z}_p} \cdots \int_{\mathbb{Z}_p} \prod_{i=0}^{n-1} (x_1 x_2 \cdots x_k-i )   d\mu_0 (x_1 )d\mu_0 (x_2 ) \cdots d\mu_0 (x_k )
\nonumber\\
&= & \int_{\mathbb{Z}_p} \int_{\mathbb{Z}_p} \cdots \int_{\mathbb{Z}_p} (x_1 x_2 \cdots x_k )_n  d\mu_0 (x_1 )d\mu_0 (x_2 )\cdots d\mu_0 (x_k )
=\check{D}_n^{(k)},
\end{eqnarray}

\noindent
we obtain the Daehee numbers of the first kind with order k, see \cite{10}.\\

\noindent
{\bf Case 6:} Setting  $x_1 x_2\cdots x_k=x$  in Eq. \eqref{12}, we obtain
\begin{equation} \label{41}
\check{D}_{n;\bar{\pmb \alpha},\bar{\pmb r}} =\int_{\mathbb{Z}_p} x-\alpha_0  )^{r_0} (x-\alpha_1  )^{r_1} \cdots (x-\alpha_{n-1}  )^{r_{n-1}} d\mu_0 (x).
\end{equation}

\noindent
{\bf Case 7:} Setting  $r_i=1$ in Eq. \eqref{37}, we obtain
\begin{equation} \label{42}
\check{D}_{n;\bar{\pmb \alpha}} =\int_{\mathbb{Z}_p} (x-\alpha_0  )(x-\alpha_1  )\cdots  (x-\alpha_{n-1}  )  d\mu_0 (x).
\end{equation}

\noindent
Which we define $\check{D}_{n; \bar{\pmb \alpha}}$  by generalized Daehee numbers of the first kind.

\begin{theorem} \label{Th:4.1}
\begin{equation} \label{43}
\int_{\mathbb{Z}_p} (x-\alpha_0  )(x-\alpha_1) \cdots  (x-\alpha_{n-1}  )  d\mu_0 (x) = \sum_{i=0}^n S(n,i; \bar{\pmb \alpha} )  \frac{s(n+i,i)}{\binom{n+i}{i}}.
\end{equation}
\end{theorem}

\begin{proof}
Substituting Eq. \eqref{9} into Eq. \eqref{42}, we obtain Eq. \eqref{43}.\\
\end{proof}

\noindent
{\bf Case 8:} Setting   $\int_{\mathbb{Z}_p}  \int_{\mathbb{Z}_p} \cdots  \int_{\mathbb{Z}_p}=\int_0^{\ell_1} \int_0^{\ell_2} \cdots \int_0^{\ell_k}$ in Eq. \eqref{12}, we obtain multiparameter Poly-Cauchy numbers of the first kind  $C_{n;\bar{\pmb \alpha},\bar{\pmb r}}^{(k)}$, see Eq. \eqref{19}.

\section{ Generalized Higher Order Daehee Numbers and Polynomials of the Second Kind}
In this section, new definition of the higher order Daehee polynomials and numbers of the second kind are established. New results are derived. Also, some special cases are introduced.
\begin{definition} \label{Def:5.1}
The generalized higher order Daehee numbers of the second kind  $\widehat{\check{D}}_{n;\bar{\pmb \alpha},\bar{\pmb r}}^{(k)}$  are defined by
\begin{equation} \label{44}
\widehat{\check{D}}_{n;\bar{\pmb \alpha},\bar{\pmb r}}^{(k)}=\int_{\mathbb{Z}_p} \int_{\mathbb{Z}_p} \cdots \int_{\mathbb{Z}_p} \prod_{i=0}^{n-1} (-x_1 x_2 \cdots x_k-\alpha_i  )^{r_i} d\mu_0 (x_1 )d\mu_0 (x_2 )\cdots d\mu_0 (x_k ).
\end{equation}
\end{definition}

\begin{theorem} \label{Th:5.1}
The numbers $\widehat{\check{D}}_{n;\bar{\pmb \alpha},\bar{\pmb r}}^{(k)}$  satisfy the relation

\begin{equation} \label{45}
\widehat{\check{D}}_{n;\bar{\pmb \alpha},\bar{\pmb r}}^{(k)}= \sum_{m=0}^{|r|}  S_{\bar{\pmb \alpha}} (n,m;\bar{\pmb r}) \sum_{\ell=0}^m L(m,n)  \sum_{\ell_1=0}^m \cdots \sum_{\ell_k=0}^m \prod_{i=0}^k S(m,\ell_i )  D_{\ell_i}.
\end{equation}

\end{theorem}

\begin{proof}
From Eq. \eqref{44}, we have
\begin{equation*}
\widehat{\check{D}}_{n;\bar{\pmb \alpha},\bar{\pmb r}}^{(k)}=\int_{\mathbb{Z}_p} \int_{\mathbb{Z}_p} \cdots \int_{|mathbb{Z}_p} \sum_{m=0}^{|r|} S_{\bar{\pmb \alpha}}  (n,m;\bar{\pmb r} ) (–x_1 x_2 \cdots x_k  )_m  d\mu_0 (x_1 )\cdots d\mu_0 (x_k )
\end{equation*}

\noindent
From definition of Lah numbers
\begin{equation*}
(–x_1 x_2 \cdots x_k  )_m =\sum_{\ell=0}^m L(m,n) (x_1 x_2 \cdots x_k  )_{\ell}.
\end{equation*}

\noindent
Hence
\begin{equation} \label{46}
\widehat{\check{D}}_{n;\bar{\pmb \alpha},\bar{\pmb r}}^{(k)}=
\int_{\mathbb{Z}_p} \int_{\mathbb{Z}_p} \cdots \int_{\mathbb{Z}_p} \sum_{m=0}^{|r|} S_{\bar{\pmb \alpha}}  (n,m;\bar{\pmb r} ) \sum_{\ell=0}^m L(m,n) (x_1 x_2 \cdots x_k  )_{\ell} d\mu_0 (x_1 )\cdots d\mu_0 (x_k ).
\end{equation}

\noindent
Substituting Eq. \eqref{13} into \eqref{46}, we obtain \eqref{45}.
\end{proof}

\begin{definition} \label{Def:5.2}
The multiparameter Poly-Cauchy numbers of the second kind  $\hat{C}_{n;\bar{\pmb \alpha}, \bar{\pmb r}}^{(k)}$ are defined by
\begin{equation} \label{47}
\hat{C}_{n;\bar{\pmb \alpha}, \bar{\pmb r}}^{(k)}=\int_0^{\ell_1} \int_0^{\ell_2} \cdots \int_0^{\ell_k} \prod_{i=0}^{n-1} (-x_1 x_2\cdots x_k-\alpha_i )^{r_i}  dx_1 dx_2 \cdots dx_k.
\end{equation}
\end{definition}

\begin{theorem} \label{Th:5.2}
The numbers $\hat{C}_{n;\bar{\pmb \alpha},\bar{\pmb r}}^{(k)}$ satisfy the relation
\begin{equation} \label{48}
\hat{C}_{n;\bar{\pmb \alpha}, \bar{\pmb r}}^{(k)}= \sum_{i=0}^{|r|} \sum_{\ell=0}^m s_{\bar{\pmb \alpha}}  (n,i;\bar{\pmb r} )L(m,n)  \frac{ (\ell_1 \ell_2 \cdots \ell_k )^{j+1}}{(j+1)^k} .
\end{equation}
\end{theorem}

\begin{proof}
From Eq. \eqref{44}, we have
\begin{eqnarray*}
\hat{C}_{n; \bar{\pmb \alpha}, \bar{\pmb r}}^{(k)}&=&
\int_0^{\ell_1}  \int_0^{\ell_2} \cdots  \int_0^{\ell_k} \sum_{m=0}^{|r|}  s_{\bar{\pmb \alpha}}  (n,m;\bar{\pmb r}) (-x_1 x_2 \cdots x_k  )^m  dx_1 dx_2 \cdots dx_k
\\
&= & \int_0^{\ell_1} \int_0^{\ell_2} \cdots \int_0^{\ell_k} \sum_{m=0}^{|r|}  s_{\bar{\pmb \alpha}}(n,m;\bar{\pmb r} ) \sum_{\ell=0}^m L(m,n) (x_1 x_2 \cdots x_k  )_\ell dx_1 dx_2 \cdots dx_k.
\end{eqnarray*}

\noindent
Then, we obtain \eqref{48}.

\end{proof}


\begin{definition} \label{Def:6.1}
The generalized higher order Daehee polynomials of the second kind $\widehat{\check{D}}_{n;\bar{\pmb \alpha},\bar{\pmb r}}^{(k)}(x)$ are defined by
\begin{equation} \label{49}
\widehat{\check{D}}_{n;\bar{\pmb \alpha},\bar{\pmb r}}^{(k)}(x)=\int_{\mathbb{Z}_p} \int_{\mathbb{Z}_p} \cdots \int_{\mathbb{Z}_p} \prod_{i=0}^{n-1} (-x_1 x_2 \cdots x_k x-\alpha_i  )^{r_i }   d\mu_0 (x_1 )d\mu_0 (x_2 ) \cdots d\mu_0 (x_k ).
\end{equation}

\end{definition}

\begin{theorem} \label{6.1}
The numbers $\widehat{\check{D}}_{n;\bar{\pmb \alpha},\bar{\pmb r}}^{(k)}(x)$ satisfy the relation
\begin{equation} \label{50}
 \widehat{\check{D}}_{n;\bar{\pmb \alpha},\bar{\pmb r}}^{(k)}(x)=
 \sum_{m=0}^{|r|} S_{\bar{\pmb \alpha}} \sum_{\ell=0}^m L(m,n)  \sum_{\ell_1=0}^m \cdots \sum_{\ell_k=0}^m \prod_{i=0}^k S(m,\ell_i ) \check{D}_{\ell_i} (x).
 \end{equation}
 \end{theorem}

\begin{proof}
From Eq. \eqref{49}, we have
\begin{eqnarray*}
\widehat{\check{D}}_{n;\bar{\pmb \alpha},\bar{\pmb r}}^{(k)}(x) &=&
\int_{\mathbb{Z}_p}  \int_{\mathbb{Z}_p} \cdots \int_{\mathbb{Z}_p} \sum_{m=0}^{|r|}  S_{\bar{\pmb \alpha}}(n,m;\bar{\pmb r})  (–x_1 x_2 \cdots x_k  )_m  d\mu_0 (x_1 ) \cdots d\mu_0 (x_k )
\\
&=& \int_{\mathbb{Z}_p}  \int_{\mathbb{Z}_p}\cdots  \int_{\mathbb{Z}_p} \sum_{m=0}^{|r|}  S_{\bar{\pmb \alpha}}(n,m;\bar{\pmb r} )  \sum_{\ell=0}^m L(m,n) (x_1 x_2 \cdots x_k  )_\ell d\mu_0 (x_1 ) \cdots d\mu_0 (x_k ).
\end{eqnarray*}

\noindent
By using Theorem \ref{Th:5.1}, the proof is completed.
\end{proof}

\begin{definition} \label{Def:6.2}
The generalized higher order Dahee polynomial  $\widehat{\check{D}}_n^{(k)}(x)$ are defined by
\begin{equation} \label{51}
\widehat{\check{D}}_n^{(k)}(x)=\int_{\mathbb{Z}_p}  \int_{\mathbb{Z}_p} \cdots \int_{\mathbb{Z}_p} (-x_1 x_2 \cdots x_k  x)_n  d\mu_0 (x_1 )\cdots d\mu_0 (x_k ).
\end{equation}
\end{definition}

\begin{theorem} \label{6.2}
The numbers
$\widehat{\check{D}}_{n;\bar{\pmb \alpha},\bar{\pmb r}}^{(k)}(x)$ satisfy the relation
\begin{equation} \label{52}
\widehat{\check{D}}_{n;\bar{\pmb \alpha},\bar{\pmb r}}^{(k)}(x)= \sum_{\ell=0}^{|r|} (-1)^\ell S(n,\ell;\bar{\pmb \alpha}, \bar{\pmb r})  \check{D}_{\ell}^{(k)} (x),
\end{equation}

\noindent
where $(x;\bar{\pmb \alpha}, \bar{\pmb r})_n=\prod_{i=0}^{n-1}(x-\alpha_i )^{r_i }$, $\bar{\pmb \alpha}=(\alpha_0,\alpha_1,\cdots, \alpha_{n-1})$, $\bar{\pmb r}=(r_0,r_1,\cdots,r_{n-1} )$.

\end{theorem}

\begin{proof}
From Eq. \eqref{49} and using Eq. \eqref{14} and Eq. \eqref{16}, we have
\begin{equation*}
\widehat{\check{D}}_{n;\bar{\pmb \alpha},\bar{\pmb r}}^{(k)}(x)=
=\int_{\mathbb{Z}_p} \int_{\mathbb{Z}_p} \cdots \int_{\mathbb{Z}_p} \sum_{m=0}^{|r|}  (-1)^m s_{\bar{\pmb \alpha}} (n,m;\bar{\pmb r} ) (x_1 x_2 \cdots x_k  x)^m  d\mu_0 (x_1 ) \cdots d\mu_0 (x_k ).
\end{equation*}

\noindent
By using Theorem \ref{Th:5.2} the proof is completed.
\end{proof}

\subsection{ Some special cases}

\noindent
{\bf Case 1:} (i) Setting  $r_i=r,\alpha_i=i$  in Eq. \eqref{49}, we have
\begin{eqnarray} \label{53}
\widehat{\check{D}}_{n;r,\bar{\pmb r}}^{(k)} (x) &= &
\int_{\mathbb{Z}_p} \int_{\mathbb{Z}_p} \cdots \int_{\mathbb{Z}_p} \prod_{i=0}^{n-1} (-x_1 x_2 \cdots x_k x-i )^r  d\mu_0 (x_1 )d\mu_0 (x_2 ) \cdots d\mu_0 (x_k ) \nonumber\\
&=& \int_{\mathbb{Z}_p} \int_{\mathbb{Z}_p} \cdots \int_{\mathbb{Z}_p} (-x_1 x_2 \cdots x_k x )_{nr}  d\mu_0 (x_1 ) d\mu_0 (x_2 ) \cdots d\mu_0 (x_k ).
\end{eqnarray}

\noindent
Replacing  $nr$ by $n$, we obtain the generalized Daehee polynomials of the second kind.\\

\noindent
(ii) Setting  $r_i=r,\alpha_i=i$, in Eq. \eqref{44}, we obtain
\begin{equation} \label{54}
\widehat{\check{D}}_{n;i,\bar{\pmb r}}^{(k)}=\int_{\mathbb{Z}_p} \int_{\mathbb{Z}_p} \cdots \int_{\mathbb{Z}_p} (-x_1 x_2 \cdots x_k  )_{nr} d\mu_0 (x_1 )d\mu_0 (x_2 )\cdots d\mu_0 (x_k ).
\end{equation}

\noindent
Replacing  $nr$ by $n$, we obtain the higher-order Daehee numbers of the second kind.\\

\noindent
{\bf Case 2:} (i) Setting  $r_i=r,\alpha_i=\alpha$ in Eq. \eqref{49}, we obtain
\begin{eqnarray} \label{55}
\widehat{\check{D}}_{n;\alpha,r}^{(k)} (x) &= &
\int_{\mathbb{Z}_p} \int_{\mathbb{Z}_p} \cdots \int_{\mathbb{Z}_p}(-x_1 x_2\cdots x_k x-\alpha )^{nr} d\mu_0 (x_1 )d\mu_0 (x_2 )\cdots d\mu_0 (x_k )
\nonumber\\
&=& \int_{\mathbb{Z}_p}  \int_{\mathbb{Z}_p} \cdots \int_{\mathbb{Z}_p} \sum_{\ell=0}^{nr} S(nr,\ell)  (-x_1 x_2 \cdots x_k x-\alpha )_{\ell}  d\mu_0 (x_1 )d\mu_0 (x_2 ) \cdots d\mu_0 (x_k )
\nonumber\\
&=& \sum_{\ell=0}^{nr} S(nr,\ell)  \int_{\mathbb{Z}_p} \int_{\mathbb{Z}_p} \cdots \int_{\mathbb{Z}_p} (-x_1 x_2 \cdots x_k x-\alpha )_{\ell}  d\mu_0 (x_1 )d\mu_0 (x_2 ) \cdots d\mu_0 (x_k )
\nonumber\\
&= & \sum_{\ell=0}^{nr} S(nr,\ell)  \widehat{\check{D}}_{\ell,\alpha}^{(k)} (x).
\end{eqnarray}

\noindent
(ii) Setting  $r_i=r,\alpha_i=\alpha$ in Eq. \eqref{44}, we obtain
\begin{eqnarray} \label{56}
\widehat{\check{D}}_{n;\alpha,r}^{(k)} &= & \int_{\mathbb{Z}_p}  \int_{\mathbb{Z}_p} \cdots \int_{\mathbb{Z}_p} (-x_1 x_2 \cdots x_k-\alpha )^{nr}  d\mu_0 (x_1 )d\mu_0 (x_2 )\cdots d\mu_0 (x_k )
\nonumber\\
&=& \sum_{\ell=0}^{nr} S(nr,\ell) \widehat{\check{D}}_{\ell,\alpha}^{(k) }.
\end{eqnarray}

\noindent
At $\alpha_i=0$ in Eq. \eqref{55}, we obtain
\begin{equation} \label{57}
\widehat{\check{D}}_{n;0,r}^{(k)} (x)=\sum_{\ell=0}^{nr} S(nr,\ell) \widehat{\check{D}}_{\ell}^{(k) } (x).
\end{equation}

\noindent
At $\alpha_i=0$ in Eq. \eqref{57}, we obtain
\begin{equation} \label{58}
\widehat{\check{D}}_{n;0,r}^{(k)}=\sum_{\ell=0}^{nr} S(nr,\ell) \widehat{\check{D}}_{\ell}^{(k)}.
\end{equation}

\noindent
{\bf Case 3:} (i) Setting  $r_i=1,\alpha_i=\alpha$ in Eq. \eqref{49}, we obtain
\begin{equation} \label{59}
\widehat{\check{D}}_{n;\alpha,1}^{(k)} (x)=\sum_{\ell=0}^n S(n,\ell) \widehat{\check{D}}_{\ell,\alpha} (x).
\end{equation}

\noindent
(ii) Setting  $r_i=1,\alpha_i=\alpha$ in Eq. \eqref{44}, we obtain
\begin{equation} \label{60}
\widehat{\check{D}}_{n;\alpha,1}^{(k)}=\sum_{\ell=0}^n S(n,\ell) \widehat{\check{D}}_{\ell,\alpha}^{(k)} .
\end{equation}

\noindent
(iii) Setting  $r_i=1,  \alpha_i=1$ in Eq. \eqref{49}, we obtain
\begin{eqnarray} \label{61}
\widehat{\check{D}}_{n;1,1}^{(k) } (x) &=&
\int_{\mathbb{Z}_p} \int_{\mathbb{Z}_p} \cdots \int_{\mathbb{Z}_p} (-x_1 x_2 \cdots x_k x-1 )^n  d\mu_0 (x_1 )d\mu_0 (x_2 )\cdots d\mu_0 (x_k )
\nonumber\\
&=& \sum_{\ell=0}^n S(n,\ell) \widehat{\check{D}}_{\ell,1}^{(k) } (x).
\end{eqnarray}

\noindent
(iv) Setting  $r_i=1,   \alpha_i=1$ in Eq. \eqref{44}, we obtain
\begin{eqnarray} \label{62}
\widehat{\check{D}}_{n;1,1}^{(k) } & =&
\int_{\mathbb{Z}_p} \int_{\mathbb{Z}_p} \cdots \int_{\mathbb{Z}_p} (-x_1 x_2 \cdots x_k-1 )^n  d\mu_0 (x_1 )d\mu_0 (x_2 ) \cdots d\mu_0 (x_k )
\nonumber\\
&=& \sum_{\ell=0}^n S(n,\ell) \widehat{\check{D}}_{\ell,1}^{(k) }.
\end{eqnarray}

\noindent
{\bf Case 4:} (i) Setting  $r_i=1,\alpha_i=0$ in Eq. \eqref{49}, we obtain
\begin{eqnarray} \label{63}
\widehat{\check{D}}_{n;0,1}^{(k) } (x) & = &
\int_{\mathbb{Z}_p}  \int_{\mathbb{Z}_p} \cdots \int_{\mathbb{Z}_p} (-x_1 x_2 \cdots x_k x)^n   d\mu_0 (x_1 ) d\mu_0 (x_2 ) \cdots d\mu_0(x_k )
\nonumber\\
&= & \sum_{\ell=0}^n S(n,\ell) \widehat{\check{D}}_{\ell }^{(k) } (x).
\end{eqnarray}

\noindent
(ii) Setting  $r_i=1,\alpha_i=0$ in Eq. \eqref{44}, we obtain

\begin{equation} \label{64}
\widehat{\check{D}}_{n;0,1}^{(k)}= \sum_{\ell=0}^n S(n,\ell) \check{D}_{\ell}^{(k)}.
\end{equation}

\noindent
{\bf Case 5:} (i) Setting  $r_i=1, \alpha_i=i,i=0,1,\cdots ,n-1$  in Eq. \eqref{61}, we have
\begin{eqnarray} \label{65}
\widehat{\check{D}}_{n;i,1}^{(k)} (x) & = &
\int_{\mathbb{Z}_p}  \int_{\mathbb{Z}_p} \cdots \int_{\mathbb{Z}_p} \prod_{i=0}^{n-1} (-x_1 x_2 \cdots x_k x-i )    d\mu_0 (x_1 ) d\mu_0 (x_2 ) \cdots d\mu_0 (x_k )
\nonumber\\
&= & \int_{\mathbb{Z}_p} \int_{\mathbb{Z}_p} \cdots  \int_{\mathbb{Z}_p} (-x_1 x_2 \cdots x_k x)_n  d\mu_0 (x_1 ) d\mu_0 (x_2 ) \cdots d\mu_0 (x_k )=\widehat{\check{D}}_n^{(k)} (x),
\end{eqnarray}

\noindent
we obtain the higher-order Daehee polynomials which defined by Kim see \cite{10}.\\

\noindent
(ii) Setting  $r_i=1, \alpha_i=i, i=0,1,\cdots,n-1$  in Eq. \eqref{44}, we have
\begin{eqnarray} \label{66}
\widehat{\check{D}}_{n;i,1}^{(k)} & = &
\int_{\mathbb{Z}_p} \int_{\mathbb{Z}_p} \cdots  \int_{\mathbb{Z}_p} \prod_{i=0}^{n-1} (-x_1 x_2 \cdots x_k-i )   d\mu_0 (x_1 ) d\mu_0 (x_2 ) \cdots d\mu_0 (x_k )
\nonumber\\
&= & \int_{\mathbb{Z}_p} \int_{\mathbb{Z}_p} \cdots  \int_{\mathbb{Z}_p} (-x_1 x_2 \cdots x_k )_n  d\mu_0 (x_1 ) d\mu_0 (x_2 ) \cdots d\mu_0 (x_k ) = \widehat{\check{D}}_n^{(k)},
\end{eqnarray}

\noindent
we obtain the Daehee numbers of the second  kind with order $k$, see \cite{10}.\\

\noindent
{\bf Case 6:} Setting  $x_1 x_2 \cdots x_k=x$  in Eq. \eqref{44}, we obtain
\begin{equation} \label{67}
\widehat{\check{D}}_{n; \bar{\pmb \alpha}, \bar{\pmb r}} =\int_{\mathbb{Z}_p} -(x-\alpha_0  )^{r_0} (x-\alpha_1 )^{r_1} \cdots  (x-\alpha_{n-1}  )^{r_{n-1}}  d\mu_0 (x).
\end{equation}

\noindent
{\bf Case 7:} Setting  $r_i=1$ in Eq. \eqref{63}, we obtain
\begin{equation} \label{68}
\widehat{\check{D}}_{n;\bar{\pmb \alpha}}=\int_{\mathbb{Z}_p} -(x-\alpha_0  )(x-\alpha_1  )\cdots  (x-\alpha_{n-1}  )  d\mu_0 (x).
\end{equation}

\noindent
Which we define $\widehat{\check{D}}_{n;\bar{\pmb \alpha}}$  by generalized Daehee numbers of the second kind.\\

\noindent
{\bf Case 8:} Setting   $\int_{\mathbb{Z}_p} \int_{\mathbb{Z}_p} \cdots \int_{\mathbb{Z}_p}=\int_0^{\ell_1}  \int_0^{\ell_2} \cdots \int_0^{\ell_k}$  in Eq. \eqref{44}, we obtain multiparameter Poly-Cauchy numbers of the second kind  $C_{n;\bar{\pmb \alpha}, \bar{\pmb r}}^{(k) }$, see Eq. \eqref{47}.


\begin{thebibliography}{99}
\bibitem{1}
Carlitz L., A note on Bernoulli and Euler polynomials of the second kind, Scripta Math. \textbf{25} (1961), pp. 323--330.

\bibitem{2}
Comtet, L., Advanced Combinatorics, Reidel, Dordrecht, 1974.
	
\bibitem{3}
Dolgy D.V., Kim T., Lee B. and Lee S.-H., Some new identities on the twisted Bernoulli and Euler polynomials, J. Comput.  Anal. Appl. \textbf{14} (2013), No. 3, pp. 441--451.
	
\bibitem{4}
El-Desouky B.S., The multiparameter non-central Stirling numbers, Fibonacci Quart. \textbf{32} (1994), No. 3, pp. 218--225
	
	
\bibitem{6}
El-Desouky B.S. and Gomaa R. S., Multiparameter Poly-Cauchy and Poly-Bernoulli numbers and polynomials, arXive:1410.5300v1 [Math. CO] 20 Oct.
	
	
\bibitem{8}
El-Desouky, B. S. and Mustafa, A., New Results and Matrix Representation for Daehee and Bernoulli Numbers and Polynomials, Applied Mathematical Sciences,  \textbf{9} (2015), No. 73, pp. 3593 -- 3610.
	
\bibitem{9}
El-Desouky B. S. and Mustafa A., New Results on Higher-Order Daehee and Bernoulli Numbers and Polynomials, Advances in Difference Equations, (2016) 2016:32. DOI 10.1186/s13662-016-0764-z
	
\bibitem{10}
Kim D.S., Kim, T., Lee S.-H. and Seo J-J. Higher-Order Daehee Numbers and Polynomials, Int. J. Math. Anal. (HIKARI Ltd), \textbf{8} (2014), No. 6, pp. 273--283.
	
\bibitem{11}
Kim D.S., Kim T., Lee S.-H. and Seo J-J.  A Note on the lambda Daehee Polynomials, Int. J.Math. Anal. (HIKARI Ltd), \textbf{7} (2013), No. 62, pp. 306--3080.
	
\bibitem{12}
Kim D.S. , Kim T. , Lee S.-H. and Seo J-J. A Note on Twisted  $\lambda$- Daehee Polynomials, Appl.Math. Sci. (HIKARI Ltd), \textbf{7} (2013), No. 141, pp. 7005--7014.
	
\bibitem{13}
Kim D.S. and Kim T., Daehee Numbers and polynomials, Appl. Math. Sci. (HIKARI Ltd.), \textbf{7} (2013), No. 120, pp. 5969--5976.
	
\bibitem{14}
Kim T. and Simsek Y., Analytic continuation of the multiple Daehee q-l-functions associated with Daehee numbers, Russ. J. Math. Phys. \textbf{15} (2008), pp. 58--65.
	
\bibitem{15}
Kimura N., On universal higher order Bernoulli numbers and polynomials, Report of the Research Institute of Industrial Technology, Nihon University, Number 70 (2003): ISSN 0386--1678.
	
\bibitem{16}
Liu G.-D. and Srivastava H. M. Explicit Formulas for the N\"{o}rlund Polynomials $B_n(x)$ and $b_n (x)$ , Computers and Mathematics with Applicatios, \textbf{51} (2006), pp. 1377--1384.
	
\bibitem{17}
Ozden H. and Cangul N. and Simsek Y., Remarks on q-Bernoulli numbers associated with Daehee numbers, Adv. Stud. Contemp. Math. (Kyungshang), \textbf{18} (2009), pp. 41--48.
	
\bibitem{18}
Wang W., Generalized higher order Bernoulli number pairs and generalized Stirling number pairs, J. Math. Anal. Appl., \textbf{364} (2010), pp. 255--274.
\end{thebibliography}
\end{document}